\documentclass[oneside,a4paper]{article}
\usepackage{pp}
\usepackage{pictexwd,dcpic}
\usepackage{multirow}
\usepackage[russian,english]{babel}
\usepackage[cp1251]{inputenc}

\def\3{\!3}
\def\4{\!4}
\def\Isom{\mathop{\rm Isom}\nolimits}

\begin{document}

\title{Telescopic actions}
\author{D. Panov\thanks{is a Royal Society University Research
Fellow}
\ and A. Petrunin\thanks{was partially supported by NSF grant DMS 0905138}}
\date{}
\maketitle

\section{Introduction}

In this paper we construct some examples of \emph{telescopic actions}
defined as follows:

\begin{thm}{Definition}\label{def:telescopic}
A co-compact properly discontinuous isometric group action $\Gamma\acts X$
on a metric space $X$
is called \emph{telescopic} if
given a finitely presented group $G$, there exists
a subgroup $\Gamma'$ of finite index in $\Gamma$ such that $G$ is isomorphic to the fundamental group of $X/\Gamma'$.
\end{thm}

Here is the first example.

\begin{thm}{Theorem}\label{thm:2D}
There is a telescopic action $\Gamma\acts X$ on a  $2$-dimensional $\CAT[-1]$ space $X$ glued from hyperbolic triangles.
\end{thm}

Denote by $\Tor \Gamma'$ the set of elements of finite order in $\Gamma'$ and
by $\<\Tor \Gamma'\>$ the subgroup of $\Gamma'$ generated by $\Tor \Gamma'$.
If $X$ is a $\CAT[-1]$ space
then given $\gamma\in\Gamma'$,
we have $\gamma\in \Tor \Gamma'$ if and only if $\gamma$ has a fixed point if $X$.
It follows that the fundamental group of $X/\Gamma'$
is isomorphic to the quotient group $\Gamma'/\<\Tor \Gamma'\>$, see \cite{armstrong}.
Therefore Theorem~\ref{thm:2D} implies  the following.

\begin{thm}{Theorem}\label{thm:group}
There exists a finitely presented hyperbolic group $\Gamma$ such that
for any finitely presented group $G$ one can find  a finite index subgroup $\Gamma'$ in $\Gamma$ such that $G$ is isomorphic to $\Gamma'/\<\Tor \Gamma'\>$.
\end{thm}

The following theorem states the existence of a telescopic action on $\HH^3$ (the 3-dimensional hyperbolic space)  with some additional properties.

Denote by $\Gamma_{\!12}$ the Coxeter group generated
by reflection in faces of a right-angled hyperbolic dodecahedron
and let $\Gamma_{\!12}\z\acts \HH^3$ be the
corresponding action.

\begin{thm}{Theorem}\label{thm:orbifold-3}
Given a finitely presented group $G$
there is a finite index subgroup $\Gamma'\subset \Gamma_{\!12}$
such that the fundamental group of $\HH^3/\Gamma'$ is isomorphic to $G$.

Moreover, the subgroup $\Gamma'\subset \Gamma_{\!12}$
can be chosen  so that
the quotient space $\HH^3/\Gamma'$ is a pseudomanifold with no boundary.
In other words, the singular points of $\HH^3/\Gamma'$
are modeled on the orientation preserving actions of $\ZZ_2$ and $\ZZ_2\oplus\ZZ_2$,
and on the action of $\ZZ_2$ by central symmetry.

\end{thm}

Note that the only topological singularities of $\HH^3/\Gamma'$
 in Theorem \ref{thm:orbifold-3} are cones over $\RP^2$ (these correspond to the centrally symmetric action of $\ZZ_2$).
That implies in particular the following result which was announced earlier by Aitchison.

\begin{thm}{Corollary}
Any finitely presented group $G$
is isomorphic to the fundamental group of $M/\ZZ_2$,
where $M$ is a closed
 oriented 3-dimensional manifold
and the action $\ZZ_2\acts M$ has only isolated fixed points.
\end{thm}

The above statement might look surprising since the fundamental groups of 3-dimensional manifolds
satisfy various severe restrictions.
For example,
\begin{itemize}
\item By a result of Heil \cite{heil}, for any $|m|\not=|n|$
the Baumslag--Solitar group  $\langle\, x,y\mid x^n\cdot y=y\cdot x^m\,\rangle$
cannot appear as a subgroup of the fundamental group of a 3-manifold.
\item
For the fundamental groups of closed 3-manifolds,
there exist algorithms to solve word problem, conjugacy problem and isomorphism problem; see the blog post of Wilton \cite{wilton} and the references therein.
\end{itemize}

\medskip

For our next result we use so-called {\it right-angled hyperbolic 120-cell}
which is a regular polytope with $120$ faces that are right-angled hyperbolic
dodecahedra (see \cite{Cox}).
Let $\Gamma_{\!120}$ be the Coxeter group generated by reflections in the faces of
the polytope and consider the corresponding action  $\Gamma_{\!120}\acts \HH^4$.

\begin{thm}{Theorem}\label{thm:orbifold}
Given a finitely presented group $G$
there is a finite index subgroup $\Gamma'\subset \Gamma_{\!120}$
such that the fundamental group of $\HH^4/\Gamma'$ is isomorphic to $G$.

Moreover the subgroup $\Gamma'\subset \Gamma_{\!120}$ can be chosen
in the index two subgroup of $\Gamma_{\!120}$ of orientation preserving
transformations.
\end{thm}

Similarly to Theorem~\ref{thm:orbifold-3}, the only topological singularities of $\HH^4/\Gamma'$ are modeled on the cone over $\RP^3$.

We use Theorem~\ref{thm:orbifold} to give an alternative short proof of the following theorem:

\begin{thm}{Taubes' theorem, \cite{T}}
\label{thm:taubes}
 For every finitely presented group $G$ there
exists a smooth compact complex $3$-manifold $W^3$ such that $\pi_1(W^3)=G$.
\end{thm}

In the original proof,
Taubes starts with an arbitrary oriented Riemannian $4$-manifold $M$ and
constructs a natural metric
on a connected sum of $M$ with sufficiently many
copies of $\overline {\CP}{}^2$.
Then he deforms the obtained metric to a metric with vanishing self-dual Weyl curvature.
This condition on the curvature tensor implies
via the Penrose construction (see \cite[13.46]{besse}), that
 the \emph{twistor bundle} over $M$
 carries a natural complex structure.
(Recall that the twistor bundle over an oriented 4-dimensional Riemannian manifold $M$
is an $\mathbb{S}^2$-bundle
with fiber over $p\in M$ formed by all isometries $J$ of the tangent space at $p$ such that $J^2=-\id$ and for which the complex orientation agrees with the given one.)
The deformation described above is the hardest part in the Taubes' proof.

We propose the following deformation-free construction.
Take a hyperbolic $4$-orbifold provided by Theorem~\ref{thm:orbifold}.
Passing to its twistor bundle we obtain a complex orbifold.
By resolving its singularities, we obtain a complex manifold with the same fundamental group.

Our proof is close in spirit to the proof of Kapovich in \cite{K},
where he shows that for any closed smooth spin $4$-manifold $M$ there exists a closed smooth $4$-manifold
$N$ such that the connected sum $M\#N$ admits a conformally flat Riemannian metric.
This implies, again by the twistor construction, that every finitely presented group is a
subgroup (in fact a free factor) of the fundamental group of a compact complex $3$-fold.

\parbf{Remarks.}
Fundamental groups of K\"ahler manifolds satisfy various non-trivial restrictions, see for example \cite{ABCKT}, and not  surprisingly
all complex manifolds obtained by our construction are non-K\"ahler (see Remark at the end of Section \ref{sec:taubes}).
In a similar vein
$3$ and $4$-dimensional hyperbolic orbifolds were used by Fine and the first author in \cite{FP}
in order to obtain non-K\"ahler manifolds with trivial canonical bundle.
We note finally, that for every finitely-presented group $G$ there exists a $2$-dimensional
irreducible complex-projective variety $W$ with the fundamental group
$G$, so that all singularities of $W$ are normal crossings and Whitney umbrellas. This was proven
very recently by Kapovich in \cite{K1} using a variation of our Theorem~\ref{thm:orbifold-3}.

\parbf{Outline of the proof.}
The first telescopic action is constructed in Section~\ref{sec:orbihedron}.

In this construction, the quotient space $Y=X/\Gamma$
 is homeomorphic to the figure eight with four attached discs;
if $g$ and $r$ are the standard generators of the figure eight, we attach the discs along the following four words: $g$,
$r$, $g{*}r$, and $g{*}r^{-1}$.
The metric inside of each disc is locally isometric to the hyperbolic plane $\HH^2$ apart from
$3$ conical points, each modeled on the singularity $\HH^2/\ZZ_2$,
and the disk boundary has geodesic curvature identically equal to $0$.

The space $X$ is constructed as the universal orbi-cover of $Y$ that has double branching points
at each of $12$ singular points;
the group $\Gamma$ is the group of deck transformations of the branched cover $X\to Y$.
Next we realize any finitely-presented group $G$
as the fundamental group of a $2$-dimensional CW-complex $Y'$,
admitting a cover $Y'\to Y$ that is allowed to double-branch
only over the $12$ singular points of $Y$.
It follows that $Y'=X/\Gamma'$ for some subgroup $\Gamma'$ of $\Gamma$.
This way we show that $\Gamma\acts X$ satisfies the telescopic property.

The actions on $\HH^3$ and $\HH^4$ are constructed in Section~\ref{sec:orbifold-3+4}.

In these constructions we use  $Y$ as a skeleton
and build orbifolds
from regular right-angled dodecahedra in $\HH^3$ and correspondingly $120$-cells from $\HH^4$.
We build them in such a way that the obtained spaces $\mathcal O_3$ and $\mathcal O_4$ naturally have a structure of a
hyperbolic orbifold and their orbi-fundamental groups admit a natural homomorphism onto $\Gamma$ with some extra properties.
Thus the universal orbi-cover of $\mathcal O_i$ is $\HH^i$
and the \emph{extra properties} ensure that the group of deck transformations has the telescopic property.

Recently, a similar construction was used by Gaifullin in \cite{Gaifullin};
he glued a compact hyperbolic $4$-manifold $M$ from hyperbolic right-angled $120$-cells such that
for any oriented compact $4$-manifold $N$ there exists a finite non-ramified cover $\tilde M$  of $M$
that admits a map of positive degree to $N$, $\tilde M\to N$.

\parbf{Acknowledgments.} In the first place we would like to thank
Joel Fine, for teaching us twistor theory, numerous discussions, and support.
We want to thank
Najmuddin Fakhruddin for the reference \cite{Ko}.
We also want to thank
Ian Agol,
Iain Aitchison,
Misha Belolipetsky
Frederic Campana,
Sergei Galkin,
Aleksandr Kolpakov,
Bruce Kleiner,
Dave Morris
 and
Ernest Vinberg
for useful conversations,
Richard Kent with Henry Wilton
for a prompt answer to our question about fundamental groups of $3$-manifolds
and Ivan Cheltsov for showing us how to resolve our singularities explicitly.
We want to thank Misha Kapovich for his interest,
thoughtful reading,
valuable suggestions,
and correcting several mistakes in the manuscript.

\section{Motivation}

The motivation comes from the following question of Gromov,
(see \cite[page 12]{gromov}).

\begin{thm}{Question}
Is it true that every compact smooth $m$-dimensional manifold $M$
is PL-homeomorphic
to the underlying space of a hyperbolic orbifold?

In other words, is there a discrete co-compact isometric action on the hyperbolic $m$-space with the quotient space PL-homeomorphic to $M$?
\end{thm}

%
%
\parit{Lower dimensions.}
It is easy to see that by passing to orbicovers of hyperbolic triangle with angles $\tfrac\pi2$, $\tfrac\pi3$ and $\tfrac\pi{5{\cdot}6{\cdot}7}$ one can get any surface.

For orientable $3$-manifolds analogous statement
is proved by Hilden, Lozano, Montesinos and Whitten in \cite{HLMW}.
They consider the hyperbolic $3$-orbifold $\mathcal{O}_3$
whose
singular locus is the Borromean rings and whose isotropy groups are all cyclic
of order four
and show that by passing to finite orbi-cover of $\mathcal{O}_3$
one can get any closed oriented 3-manifold.
The orbifold $\mathcal{O}_3$ was first considered by Thurston;
it can be obtained from the regular hyperbolic right-angled dodecahedron by gluing $6$ pairs of adjacent faces.
It seems that if instead of $\mathcal{O}_3$, one starts with the regular hyperbolic right-angled dodecahedron then one can get any (not necessary orientable) closed $3$-manifolds.

All this suggests the following variation of Gromov's question.

\begin{thm}{Question}\label{conj:telescop}
Given a positive integer $m$, is there an $m$-dimensional hyperbolic orbifold $\mathcal O_m$, such that any compact smooth $m$-dimensional manifold $M$ is  PL-homeomorphic
to the
underlying space of a finite orbi-cover of $\mathcal O_m$?

In other words,
is there a co-compact isometric discrete action on the hyperbolic $m$-space $\Gamma\acts \HH^m$,
such that $M$ is  PL-homeomorphic to $\HH^m/\Gamma'$ for some finite index subgroup $\Gamma'$ of $\Gamma$?
\end{thm}

Our construction of the actions might be considered as a solution of a further variation of this conjecture, which takes into account only the fundamental group.

\parit{Higher dimensions.}
The following question seem to be completely open.

\begin{thm}{Question}\label{quest:no-sc}
Let $m$ be a large integer.
Is there any cocompact isometric properly discontinuous action $\Gamma\acts\HH^m$
such that the quotient space $\HH^m/\Gamma$ is simply connected?

Equivalently, is there a cocompact lattice in $\Isom^+\HH^m$ generated by elements of finite order?
\end{thm}

Here $\Isom^+\HH^m$ stays for the group of orientation preserving isometries of $\HH^m$.

A negative answer would imply that there is no telescopic action on $\HH^m$ for large $m$ (because the trivial group could not be realized).

A negative answer  to Question~\ref{quest:no-sc}
would also imply
a negative answer to Gromov's question,
but much less would be sufficient.

First note the following.

\begin{thm}{Claim}
Let $\Gamma$ act isometrically and properly discontinuously on $\HH^m$ or $\RR^m$ and let $X$ be the quotient space.
Then
\begin{enumerate}
\item\label{sc} $X$ is simply connected
if and only if $\Gamma$ is generated by elements of finite order.
\item\label{PL-mnfld} If $X$ is PL-homeomorphic to a simply connected manifold
then $\Gamma$ is generated by rotations around subspaces of codimension $2$.
\end{enumerate}
\end{thm}

The part \ref{sc} follows from \cite{armstrong}.
The second part seems to be noted by Schwarzman in \cite{schwarzman}\footnote{We were not able to find this paper, but this can be proved along the same lines as \cite{armstrong}.}.
The converse for part \ref{PL-mnfld} for finite groups was proved by Mikhailova in \cite{mikhailova}.

Note that the cone over spherical suspension over  Poincar\'e sphere is homeomorphic to $\RR^5$ and it is a quotient of $\RR^5$ by a finite subgroup of $SO(5)$.
Hence, in part \ref{PL-mnfld},
one can not exchange ``PL-homeomorphism'' to ``homeomorphism''.

If the answer to Gromov's question is ``yes'',
then in particular one has to be able to construct a
hyperbolic orbifold with underlying space PL-homeomorphic to $\mathbb S^m$.
Taking above claim into account this would imply a positive answer to the following question.

\begin{thm}{Question}\label{conj:no-sc'}
Let $m$ be a large integer.
Is there a cocompact lattice in $\Isom^+\HH^m$ which is generated by
rotations around subspaces of codimension 2?
\end{thm}

Note that the orientation preserving part of any Coxeter's action,
is generated by rotations.
The non-existence of compact hyperbolic Coxeter polytopes proved by Vinberg (see \cite{V1} and \cite{V2})
suggests that the answer should be ``no''.

\section{Telescopic orbihedron.}\label{sec:orbihedron}

In this section we prove Theorem~\ref{thm:2D}.

Denote by $\star$ the $0$-cell of
the figure eight and by $g$ and $r$ its loops
($g$ is for ``green'' and $r$ is for ``red'').

First, let us construct the space $Y$ that will serve further as $X/\Gamma$.
Attach to the figure eight
four discs $\mathcal{B}$, $\mathcal{W}$, $\mathcal{G}$, $\mathcal{R}$
(named for ``black'', ``white'', ``green'', and ``red'')
along $g{*}r^{-1}$, $g{*}r$, $g$ and $r$ respectively.
It is easy to see that $Y$ is homeomorphic to $\RP^2$
with two discs attached along two lines;
$\RP^2$ is colored in black and white and the attached discs are red and green.

We equip
$Y$ with an intrinsic metric such that each disc $\mathcal{B}$, $\mathcal{W}$, $\mathcal{G}$, $\mathcal{R}$
is isometric to a disc obtained by gluing two copies of a right-angled hyperbolic pentagon along $4$ sides.
This way each disc contains three singular points
modeled on the singularity $\HH^2/\ZZ_2$.
In total we have $12$ such special points $\{p_1,p_2,\dots,p_{12}\}$
that will be the only branching points in $Y$; each $p_i$ has branching order $2$.

The space $Y$ admits the unique cover $X\to Y$ with $\CAT[-1]$
total space $X$ which has double branching at each $p_i$
(see \cite{haefliger} for details).
We let $\Gamma\acts X$ be the action of deck transformations for the cover $X\to Y$.

A cover $Y'$ of  $Y$ that might have double branching only at $p_i$ is called
an {\it orbi-cover} of $Y$.
Any such $Y'$ can be obtained as the quotient $X/\Gamma'$
for some subgroup $\Gamma'$ of $\Gamma$.
The index $[\Gamma:\Gamma']$ is the degree of the cover $Y'\to Y$.

Taking all the above into account, Theorem~\ref{thm:2D}
boils down to the following.

\begin{thm}{Proposition}\label{prop}
Given a finitely presented group $G$ there is a finite orbi-cover $f\:Y'\to Y$ such that
$\pi_1(Y')$
is isomorphic to $G$.
\end{thm}

\parit{Proof.}
For a given group $G$ we will construct a special 2-dimensional CW-complex $Y'$
with the fundamental group $G$, that admits an orbi-cover $f\:Y'\z\to Y$.
We divide the proof into two steps.
In step 1, we construct $Y'$  by attaching a finite number of discs to a closed surface.
In step 2, we construct an orbi-cover $f\:Y'\to Y$.

\parit{Step 1.}
Note that $G$ can be realized as the fundamental group of an oriented surface,
say $\Sigma_0$,
with finite number of attached discs.
Specifically, assume $G$ has $k$ generators.
Take the oriented surface $\Sigma_0$ of genus $k$.
By attaching $k$ discs to $\Sigma_0$ one can reduce its fundamental group to $\FF_k$,
the free group with $k$ generators.
Attaching further disks to $\Sigma_0$ corresponding to the relators
in $G$ one obtains a space with
the fundamental group $G$.

Let us draw on $\Sigma_0$ in red
the curves along which the discs were attached.
We may assume that all these curves intersect transversally.
We also may assume that the red curves cut $\Sigma_0$ into discs
and each curve intersects some other curve.
(That is easy to arrange by adding a finite number of null-homotopic red curves.
Attaching a disc along such a curve does not change the fundamental group of the space.)

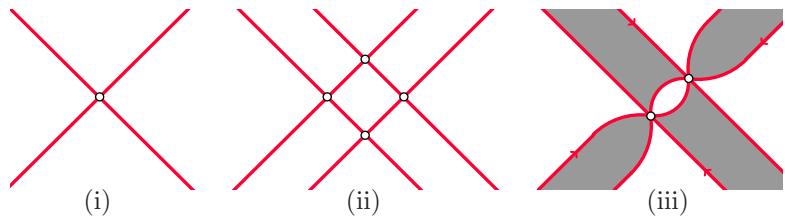
\begin{center}
\begin{lpic}[t(0mm),b(3mm),r(0mm),l(0mm)]{pics/new-42223(0.5)}
\lbl[t]{23,0;(i)}
\lbl[t]{93,0;(ii)}
\lbl[t]{173,0;(iii)}
\end{lpic}
\end{center}

For each of these curves, let us draw a parallel red curve,
so instead of one intersection as in  figure (i),
we get four intersections as in  figure (ii).
Further, deform each configuration as in figure (ii) to that in figure (iii).
Now the curves have only triple intersection points and all curves are mutually transversal.
We are about to explain the meaning of black-and-white colors and the orientation of the curves in (iii).

Note that the red curves still cut $\Sigma_0$ into discs.
Moreover now we can color the discs in the checkerboard order;
i.e., make them black and white in such a way that the disc changes its color
each time one crosses (transversally)  a red curve.
Color in black all {\it long and thin} disks whose boundary contains two pieces of parallel red curves.
Color the rest of the surface in white.
Since $\Sigma_0$ is oriented, we can orient the boundary of all black disks clockwise,
so the boundary of any white disc will be oriented counter-clockwise.

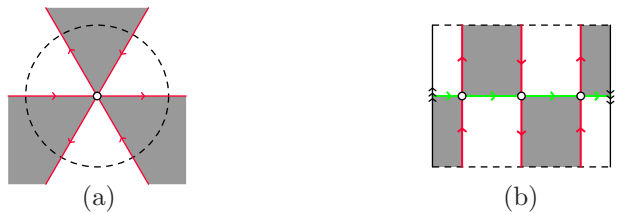
\begin{center}
\begin{lpic}[t(0mm),b(3mm),r(0mm),l(0mm)]{pics/moebius(0.5)}
\lbl[t]{24,0;(a)}
\lbl[t]{135,0;(b)}
\end{lpic}
\end{center}

Now cut from $\Sigma_0$ a small disc around each point of intersection, along the dashed line as on figure (a);
then glue instead a M\"obius band with central line marked in green, as on the figure (b).
This way we get a non-oriented surface $\Sigma$ with a net of red and green closed curves
which satisfies the following conditions:

\begin{enumerate}
\item Each curve intersects with at least one other curve
 and the intersections are transversal.
Two curves of the same color can not intersect.

\item\label{lem:X:even} The orientation on each curve can be chosen in such a way that if one goes along one of the curves then others cross it alternately from right to left and from left to right.

\item\label{lem:X:discs}
The red and green curves cut $\Sigma$ into discs.
These discs can be colored in the checkerboard order in such a way that
if one moves around the boundary of white (black) disc
then red and green segments have the same (correspondingly the opposite) orientation.

\item\label{lem:X:+discs} If one attaches a disc to each of the red and green curves then the fundamental group of the obtained space $Y'$ is isomorphic to $G$.
\end{enumerate}

Let us construct $Y'$ as it is described in condition \ref{lem:X:+discs} and
color the attached discs into green and red accordingly to the color of their boundary curve.

\parit{Step 2.}
Now let us construct a map $f\: Y'\to Y$.
Map all points of intersection of the red and green curves on $Y'$ to
$\star\in Y$
and send by one-to-one orientation preserving maps all red and green segments
of red and green curves to $r$ and $g$ correspondingly.
From property \ref{lem:X:discs},
it follows that one can extend this map to the whole $Y'$ in such a way that two-cells are mapped to the two-cells of $Y$ with the same color.
This map is homotopic to a branched cover
with branchings only over $\{p_1,p_2,\dots,p_{12}\}$ of order at most 2.
The later statement follows from the following lemma.

\begin{thm}{Lemma}\label{lem:2-points}
Let $\mathbb{D}$ be the two-dimensional disc.
Then any cover $\partial \mathbb{D}\to\partial \mathbb{D}$
can be extended to a ramified covering $\mathbb{D}\to \mathbb{D}$
which is branching only at the given
two interior points with order at most 2.
\end{thm}

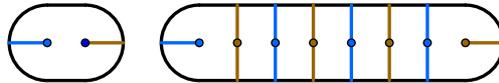
\begin{center}
\begin{lpic}[t(0mm),b(0mm),r(0mm),l(0mm)]{pics/cover(0.5)}
\end{lpic}
\end{center}

The proof should be clear from the picture; cf. \cite[Proposition 1]{feighn}. \qeds

\section{Telescopic orbifolds}\label{sec:orbifold-3+4}

We will construct the telescopic action  on $\HH^3$
and will use it further to construct the action on $\HH^4$.
Let us first give an outline
 of the construction and
then describe each case in more details.

Note that the terms ``hyperbolic orbifold''
and ``discrete isometric group action on hyperbolic space''
have the same meaning,
but in our constructions it is more intuitive to use the orbifold terminology.
(The reader has to get used to the translations from one terminology to the other;
for example,
\emph{``subaction''}
corresponds to \emph{``orbi-cover''} and so on.)

Given a hyperbolic orbifold $\mathcal O=\Gamma\curvearrowright\HH^n$,
we denote by $|\mathcal O|$ its underlying space;
i.e., $|\mathcal O|=\HH^n/\Gamma$.

The space $Y$ constructed in Section~\ref{sec:orbihedron}
will be also treated as an \emph{orbihedron};
i.e., we write $Y$ for the action $\Gamma\curvearrowright X$  and $|Y|$ for the quotient space $X/\Gamma$.
For any subgroup $\Gamma'\le\Gamma$ there is a covering map $Y'\to Y$
from $Y'=\Gamma'\curvearrowright X$ to $Y$ with branching only at the points $p_i$;
we will call $Y'$ an {\it orbi-cover} of $Y$.
(The metrics on $X$ and $Y$ constructed above will not be used further.)

To prove Theorems  \ref{thm:orbifold-3} or \ref{thm:orbifold}
we have to construct a (three or four dimensional) hyperbolic orbifold $\mathcal O$
such that every finitely presented group $G$ appears as $\pi_1|\mathcal O'|$ for a finite orbi-cover $\mathcal O'\to\mathcal O$
which satisfies the additional properties stated in the theorems.

Any orbifold $\mathcal O$ of that type will be called \emph{telescopic}.
It is straightforward to check that a hyperbolic orbifold  $\mathcal O$
is telescopic if it satisfies the following two conditions:
\begin{enumerate}
\item\label{conditions} There is an embedding $\iota\:|Y|\hookrightarrow|\mathcal O|$.
\item For any orbi-cover $Y'\to Y$
there is an orbi-cover $\mathcal O'\to O$
and an embedding
$\iota'\:|Y'|\hookrightarrow|\mathcal O'|$ such that
the following diagram is commutative
$$\raisebox{-0.9cm}{$\begindc{\commdiag}[10]
\obj(0,3)[aa]{$|Y'|$}
\obj(4,3)[bb]{$|\mathcal O'|$}
\obj(0,0)[cc]{$|Y|$}
\obj(4,0)[dd]{$|\mathcal O|$}
\mor{aa}{bb}{$\iota'$}
\mor{aa}{cc}{}
\mor{bb}{dd}{}
\mor{cc}{dd}{$\iota$}
\enddc$}
$$
and $\iota'$ induces an isomorphism
$\pi_1|Y'|\to\pi_1|\mathcal O'|$.
\end{enumerate}

The construction of $\mathcal O$ from $Y$ will be given in three steps.
To visualize the first two steps in the  construction
it is convenient to pass to a double branched cover $Y_2$ of $Y$ which we are about to describe.

\parbf{Double orbi-cover of $\bm{Y}$.}
Let us describe a double orbi-cover $Y_2$ of $Y$ that will be used further.
We realize $|Y_2|$ topologically as a cell complex in $\mathbb{S}^3$.
Namely denote by $z_{\text{red}}$ and $z_{\text{green}}$ two opposite  poles in $\mathbb{S}^3$, and let
$\mathbb{S}_{\text{black-or-white}}^2$ be the equatorial sphere.
Let $\mathbb{S}^1_{\text{red}}$, $\mathbb{S}^1_{\text{green}}$ be two great
 orthogonal circles on
$\mathbb{S}_{\text{black-or-white}}^2$.
Let $\mathbb D_{\text{red}}$ and $\mathbb D_{\text{green}}$ be the two-dimensional hemispheres in $\mathbb{S}^3$
whose centers are $z_{\text{red}}$ and $z_{\text{green}}$, and whose boundaries are $\mathbb{S}^1_{\text{red}}$ and $\mathbb{S}^1_{\text{green}}$ respectively.
With these notations
$|Y_2|$ is the union of the two disks $\mathbb D_{\text{red}}$, $\mathbb D_{\text{green}}$ and the sphere $\mathbb{S}^2_{\text{black-or-white}}$.
Finally, let $\sigma$ be the involution on $\mathbb{S}^3$ that fixes the poles
and restricts to the central symmetry on $\mathbb{S}^2_{\text{black-or-white}}$.

It is clear that the quotient $|Y_2|/\sigma$ is homeomorphic to $|Y|$.
The two intersections $\star_1$ and $\star_2$ of $\mathbb{S}^1_{\text{red}}$ with $\mathbb{S}^1_{\text{green}}$ on $|Y_2|$ correspond to the point $\star$ on $|Y|$.
The black and white two-cells of $|Y|$ correspond to $\mathbb{S}^2_{\text{black-or-white}}$,
and the red and green two-cells correspond to $\mathbb D_{\text{red}}$ and $\mathbb D_{\text{green}}$.

The one-skeleton
of $|Y_2|$ is the graph with two vertices $\star_1$ and $\star_2$ joined by
$4$ edges.
Note that $|Y_2|$ is obtained from the skeleton  by attaching $6$ two-cells;
two black,
two white,
one red and one green.
Further, $|Y_2|$ cuts from $\mathbb{S}^3$ four balls
and each two-cell of $|Y_2|$ lies in the boundary of 2 of these balls.

\parbf{Step 1: Pentagonalization.}
We glue $Y$ from pentagons in a specific way
and equip $Y$ with an intrinsic metric
such that each pentagon is isometric to a regular right-angled pentagon in $\HH^2$.
One could also think about this step as of gluing $Y_2$ from pentagons in a $\sigma$-invariant way.

The ``pentagolizations'' which we construct satisfy some additional properties that will
permit us to do the next steps in the construction;
for example, the total angle around
each of the branching points $\{p_1,p_2,\dots,p_{12}\}$ has to be equal to $\pi$.
We stress here that the metric on $Y$ induced by the pentagonalization will
differ from the one used in the proof of Theorem \ref{thm:2D}, in particular
it will have more metric singularities.
But the orbihedron structure will be identical to that in Theorem \ref{thm:2D}.

\parbf{Step 2: Attaching the meat.}
In this  step we describe a way to glue a number of \emph{hyperbolic right-angled dodecahedra}
or correspondingly
\emph{hyperbolic right-angled regular 120-cells} to the pentagons in $Y$
to obtain a telescopic orbifold with nonempty boundary.

By our construction we obtain an orbifold that corresponds to a subaction of
the action $\Gamma_{\!12}\acts \HH^3$ or the action $\Gamma_{\!120}\acts \HH^4$
correspondingly.

\parbf{Step 3: Doubling.}
In this step we get rid of the boundary by applying the doubling of the obtained orbifold across its boundary.

Again by our construction we obtain an orbifold that corresponds to a subaction of
the action $\Gamma_{\!12}\acts \HH^3$ or the action $\Gamma_{\!120}\acts \HH^4$
correspondingly.

Recall that the \emph{doubling} of a space $X$ across a
subset $A\subset X$
is obtained by gluing two copies of $X$ at the corresponding points  of the copies of $A$.
It  is easy to see that the doubling of an orbifold across its boundary carries a natural orbifold structure.

If $W$ is the doubling of $X$ across $A$ then $X$ admits two natural embeddings $l,r\:X\hookrightarrow W$,
which we call \emph{left} and \emph{right embeddings}.\label{left-embedding}

\medskip

Now we turn to the the details of the above construction in 3- and 4-dimensional cases.
You should already see the home through the woods
and it should be clear that you can get there,
we are about to describe a trail.

\subsection*{The construction of 3-orbifold}

\parbf{Pentagonalization.}
The pentagonalization we are about to construct
is different from the one in Section~\ref{sec:orbihedron}.

The pentagonalization of $\mathbb{S}^2_{\text{black-or-white}}$ in $Y_2$
is obtained by doubling of the left part of the following diagram across its boundary.
Both red and green two-cells of $Y_2$ are glued from $8$ pentagons as shown on the right diagram.
They will be attached along the corresponding
 lines on the left diagram.
The poles are marked by $z$ and the points corresponding to $\star_1$ and $\star_2$ are marked by $\star$.
The meaning of dashed lines and
blue
 and purple
points will be explained below.

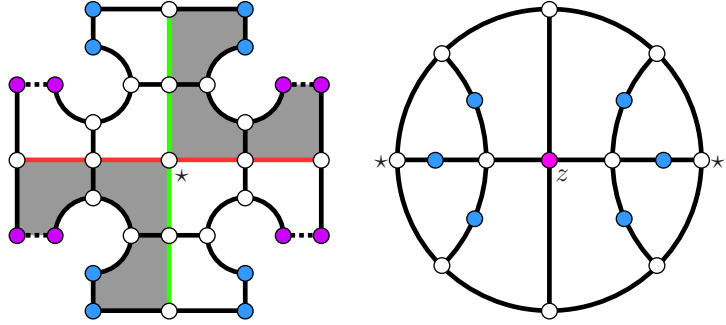
\begin{center}
\begin{lpic}[t(0mm),b(0mm),r(0mm),l(0mm)]{pics/3d-pentagolization-super-new(0.5)}
\lbl[lt]{45,42;$\star$}
\lbl[lt]{145,42;$z$}
\lbl[r]{101,44;$\star$}
\lbl[l]{186,44;$\star$}
\end{lpic}
\end{center}

In the corresponding pentagonalization of $Y$ (i.e., after taking quotient by $\sigma$),
the black and white two-cells are glued from 6 pentagons each,
and the red and green two-cells are glued from 4 pentagons each.

To specify the orbi-structure of $Y$, we need to choose three branching points $(p_i)$ on each cell out of blue
and purple
points on the diagram.
We make a choice in such a way that
on each two-cell of $Y$ one purple point is left (this point is not treated as an orbifold point of $Y$, it just represents a {\it metric singularity}).
The reason for making such a choice will become clear later on.

\parbf{Attaching the meat.}
Recall that the complement of $Y_2$ in  $\mathbb{S}^3$ is a union of four balls
and each pentagon of $Y_2$ belongs to the boundary of two balls.
We will attach dodecahedra to pentagons assigning to each dodecahedron
one of four balls to which it ``belongs''.
To each pentagon two dodecahedra will be attached
and each dodecahedron is attached to one or two pentagons.

To do so,
first consider all pairs consisting
of one green-or-red
and one black-or-white
pentagons in $Y_2$ that share one edge and belong to the boundary of the same ball.
To each such a pair we attach
a right angled hyperbolic dodecahedron
along two adjacent faces.
After that each green-or-red pentagon
and each black-or-white pentagon adjacent to the center of the left diagram
has
{\it two} dodecahedra attached.
For each remaining black-or-white pentagons we attach one dodecahedron from
the side from which it was not yet attached.

Further we glue together attached dodecahedra along pairs of faces
that intersect $Y_2$ in a common edge.
To be glued the faces of dodecahedra must have a common edge in $Y_2$
and yet satisfy
one of the following mutually exclusive conditions:

\begin{enumerate}
\item Two dodecahedra correspond to the same ball and the edge is marked by a solid line
(of any color)
  on the diagram.
\item Two dodecahedra correspond to two different balls and the edge is marked by a dashed line on the diagram.
\end{enumerate}

After these gluings, all white points and the point $z$ on the diagram become regular;
i.e., they all admit a neighborhood isometric to an open set in $\HH^3$.
The blue
points lie on a singular
 line, perpendicular to the plane of the diagram;
this line has conical angle $\pi$ around it
and therefore the corresponding singularity
 is modeled on the orientation preserving action $\ZZ_2\acts\HH^3$.
All the purple
 points, except $z$, lie at the ends of dashed lines and
they become isolated singularities
  modeled on the action $\ZZ_2\acts\HH^3$ by central symmetry.
Indeed, a simple loop on $Y_2$ encircling a
purple point represents an orientation reversing path
in the space obtained after gluing (by construction,
the normal to $Y_2$ changes its direction along such a path).

As a result, we obtain a space glued from regular right-angled dodecahedra with an isometric
 involution $\sigma$.
This space has a natural structure of hyperbolic
orbifold.
Each vertex on the diagram that is an end of a dashed line corresponds to a singularity
modeled by the action of $\ZZ_2$ by central symmetry; at the boundary of the orbifold we have orientation reversing actions of
$\ZZ_2$, $\ZZ_2\oplus\ZZ_2$ and $\ZZ_2\oplus\ZZ_2\oplus\ZZ_2$, and the rest of
the singularities are given  by orientation preserving actions of $\ZZ_2$.

Taking the quotient of this orbifold by $\sigma$, we get a new 3-dimensional orbifold,
say  $\mathcal P_3$;
it has two more singularities
at the images of
$z_{\text{red}}$ and $z_{\text{green}}$,
both modeled on the action of $\ZZ_2$ by central symmetry.
(That is why we color $z$ in purple
 on the diagram.)

By our construction
\begin{itemize}
\item $|Y|$ is sitting naturally
inside $|\mathcal P_3|$. (Note however that the inclusion $|Y|\z\hookrightarrow|\mathcal P_3|$ is not induced by a legitimate embedding in the orbihedra-category.)

\item Each branching point $p_i\in Y$ belongs to the singular locus of $\mathcal P_3$ modeled on a $\ZZ_2$ action; either by central symmetry or by reflection in a line.
The choice of $12$ points was done in such a way
that on each two-cell of $Y$ (black, white, green and red) we have one
singular point of $\mathcal P_3$ modeled on the action of $\ZZ_2$ by central symmetry
which is not a branching point of $Y$; i.e., not one of $p_i$.
\end{itemize}

In order to prove that $\mathcal P_3$ is telescopic we need to show that
each orbi-cover $Y'\to Y$ can be lifted to an orbi-cover $\mathcal P_3'\to\mathcal P_3$ such that the two
conditions on page \pageref{conditions} hold.
To construct the lifting, note that $|Y|$ is a strong deformation
retract of $|\mathcal P_3|$.
Moreover, the retraction $s\:|\mathcal P_3|\to |Y|$ can be chosen so that each preimage $s^{-1}(p_i)$ is formed by the edge(s) of dodecahedra which touch $|Y|$ at $p_i$.
(There might be two or one of such edges depending on the color
 of $p_i$.)

Then $|\mathcal P_3'|$ is obtained as follows.
Set
$$S=|\mathcal P_3|\backslash \left(\bigcup_i s^{-1}(p_i)\right).$$
Then take the \emph{fiber product} $S'=|Y'|\times_{|Y|} S$
for the map $s$
and define $|\mathcal P_3'|$ as the metric completion of $S'$.
Again, $|Y'|$ is a strong deformation retract of $|\mathcal P_3'|$.
It is clear that $|\mathcal P_3'|$ constructed this way is an underlying
space of a hyperbolic orbifold $\mathcal P_3'$ and the natural projection $|\mathcal P_3'|\to |\mathcal P_3|$ is induced by an orbi-cover $\mathcal P_3'\to \mathcal P_3$.

\parbf{Doubling.}
Let $\mathcal O_3$ be the doubling of $\mathcal P_3$ across its boundary.
From the above it follows that all singularities of $\mathcal O_3$ are modeled on the orientation preserving actions of $\ZZ_2$ and $\ZZ_2\oplus\ZZ_2$
or the action of $\ZZ_2$ by central symmetry.

Given an orbi-cover $Y'\to Y$
and the corresponding orbi-cover $\mathcal P_3'\to\mathcal P_3$
the doubling $\mathcal O_3'$ of $\mathcal P_3'$
is the total space of the orbi-cover $\mathcal O_3'\to\mathcal O_3$.

Let us show finally that the left embedding $ |\mathcal P_3'|\hookrightarrow|\mathcal O_3'|$
(defined on page \pageref{left-embedding})
induces an isomorphism
$\pi_1 |\mathcal P_3'|\z\to\pi_1 |\mathcal O_3'|$
and therefore  $\mathcal O_3$ is telescopic.

From existence of retraction of the double
$|\mathcal O_3'|$ to its left side $|\mathcal P_3'|$, it follows that
the map  $\pi_1 |\mathcal P_3'|\to\pi_1|\mathcal O_3'|$
is injective.
On the other hand, any loop in $|\mathcal O_3'|$
can be pushed inside the left copy of $|\mathcal P_3'|$ in $|\mathcal O_3'|$;
i.e. the map $\pi_1 |\mathcal P_3'|\to\pi_1|\mathcal O_3'|$ is also surjective. Let us prove this.

First, deform the  loop so that it intersects the right image of $|Y'|$ transversally in the
interiors of its two-cells.
Next, homotopy this loop further into a loop that does not intersect the right image of $|Y'|$
at all.
The later is possible by the lemma below
since each two-cell of $|Y'|\subset |\mathcal P_3'|$ has a cone point over $\RP^2$;
such a point exists, since each two-cell in $Y$ contains a cone point over $\RP^2$
which is not a branching point for the orbi-cover $Y'\to Y$.

\begin{thm}{Lemma}\label{lem:ex}
Let $\ell$ be a line in $\RP^2$.
Consider the cone $A$ over $\RP^2$ with the tip $o$
and let
and $B\subset A$ the cone over $\ell$
(again with the tip $o$).
Then any path $\gamma$ with ends $x,y\notin B$
is homotopic rel ends to a path which does not intersect $B$.
\end{thm}

\parit{Proof.}
Since $A\backslash B$ is connected, there exists a path $\gamma'$ in $A\backslash B$ connecting
$y$ to $x$.
Then the assertion of lemma is equivalent to the claim that the
concatenation of $\gamma'$ and $\gamma$
is a loop nul-homotopic in $A$.
The latter follows
from contractibility of $A$.
\qeds

One can check
explicitly that $|\mathcal P_3|\setminus |Y|$ is a product of a half-closed interval with
a  surface (in fact, a Klein bottle) and  $|\mathcal P_3'|\setminus |Y'|$ covers $|\mathcal P_3|\setminus |Y|$
with ramifications along a collection of vertical half-closed intervals.
Therefore the space  $|\mathcal P_3'|\setminus |Y'|$
is a direct product of a surface with a half-closed interval.

So, once the loop is disjoint from the right image of  $|Y'|$ we can push it inside the left copy of $|\mathcal P_3'|$. This finishes the proof of surjectivity.
\qeds

\subsection*{The construction of 4-orbifold}

We will use the action provided by Theorem~\ref{thm:orbifold-3} to construct an action required by Theorem~\ref{thm:orbifold}.
In order to resolve the ambiguity in the notation,
we denote by $\Gamma_{\3} \acts\HH^3$ and $\Gamma_{\4} \acts\HH^4$ the actions in  Theorems \ref{thm:orbifold-3} and \ref{thm:orbifold} correspondingly.

First, let us extend the action $\Gamma_{\3}\acts\HH^3$ to $\HH^4$.
Consider the hyperboloid model in $\RR^{4,1}$ for $\HH^4$.
Choose an embedding $\HH^3\z\hookrightarrow\HH^4$ which corresponds to a coordinate embedding
$\RR^{3,1}\hookrightarrow\RR^{4,1}$.
The action of $\Gamma_{\3}\acts\HH^3$ lifts to a unique representation $\Gamma_{\3}\acts\RR^{3,1}$
that does not swap two connected components of the light cone in $\RR^{3,1}$.
Denote by $A_\gamma$ the $4\times 4$ matrix which corresponds to $\gamma\in\Gamma_{\3}$.

Consider the representation $\Gamma_{\3}\acts\RR^{4,1}$ given by the block-diagonal matrix
$$\gamma\mapsto
B_\gamma
\df
\left(
\begin{array}{c|c}
\raisebox{-20pt}{{\huge\mbox{{$A_\gamma$}}}}& 0 \\ [-5.5ex]
& 0\\[-0.5ex]
& 0\\[-0.5ex]
& 0\\ [0.5ex]\hline
\raisebox{-2pt}{$0\ 0\ 0\ 0$} & \raisebox{-2pt}{$\det A_\gamma$}
\end{array}
\right)$$
Note that $\det A_\gamma=\pm1$ and $\det B_\gamma=1$;
i.e., the constructed action $\Gamma_{\3}\acts\HH^4$ is orientation preserving.

Consider the tiling of $\HH^4$ by regular right-angled $120$-cells that extends
the tiling of $\HH^3$ by dodecahedra.
Let $W$ be the union of all $120$-cells touching $\HH^3$.
Note that $W$ is an infinite Coxeter polytope.
In particular it has a natural orbifold structure.
Further, $W$ is an invariant set of the constructed action $\Gamma_{\3}\z\acts\HH^4$.

Note that $\HH^3/\Gamma'$ is a deformation retract of $W/\Gamma'$
for any subgroup $\Gamma'\le \Gamma_{\3}$.
Since $\Gamma_{\3}\acts\HH^3$ is telescopic, we get that $\mathcal P_4=\Gamma_{\3}\acts W=\Gamma_{\3\frac12}\acts \HH^4$ is telescopic; the action $\Gamma_{\3\frac12}\acts \HH^4$ is generated by elements of $\Gamma_{\3}$ and the Coxeter group of $W$.

Consider finally the orbifold $\mathcal O_4=\Gamma_{\4}\acts\HH^4$,
where $\Gamma_{\4}$ is formed by all orientation
preserving elements in $\Gamma_{3\frac12}$.

Since all orientation reversing elements of $\Gamma_{3\frac12}$ are generated by reflections in
the faces of $W$,
one can also view $\mathcal O_4$ as the doubling of $\mathcal P_4$ across the boundary;
i.e., in the subset of all points of $\mathcal P_4$  whose stabilizer includes a reflection in a hyperplane.
The left embedding $|\mathcal P_4'|\to |\mathcal O_4'|$
(defined on page \pageref{left-embedding})
induces an isomorphism $\pi_1|\mathcal P_4'|\z\to \pi_1|\mathcal O_4'|$.
The later is proven the same way as in the 3-dimensional construction, the argument is
even simpler.
It is sufficient to know that $\HH^3/\Gamma'$ is a three-dimensional orbifold
that contains at least one point of  $|\mathcal O_4'|$ modeled by the action of $\ZZ_2$ by central symmetry.

Hence $\mathcal O_4=\Gamma_{\4}\acts\HH^4$ is telescopic and the remaining conditions follow directly from the construction.

\section{Taubes' theorem}\label{sec:taubes}

Recall that the twistor space of $\mathbb{S}^4$ with
the standard conformal structure is $\CP^3$ with its standard holomorphic
structure (see \cite[13.65]{besse}).
The group $\SO(5,1)$ of conformal (orientation preserving) transformations of $\mathbb{S}^4$
acts on the twistor space by biholomorphisms, i.e., by complex
projective transformations.

For the standard conformal embedding $\HH^4\hookrightarrow \mathbb{S}^4$,
the group of conformal transformation of $\mathbb{S}^4$
preserving $\HH^4$ coincides the group of isometries of $\HH^4$. In particular
to each compact oriented hyperbolic orbifold $\HH^4/\Gamma$ corresponds
a $3$-dimensional compact complex orbifold that can be obtained by taking the quotient of the part of $\CP^3$ over $\HH^4$ by $\Gamma$; the complex orbifold
is naturally mapped to the hyperbolic one and all the fibers are topologically $\mathbb{S}^2$.

\begin{thm}{Theorem \cite[7.8.1]{Ko}} \label{kollar} Let  $V$
be a  normal  analytic  space  and  let  $f\:W\z\to V$  be  a  resolution  of
singularities.
Assume $V$ has only quotient singularities then $\pi_1(W)\z\cong \pi_1(V)$.
\end{thm}

\parit{Proof of Taubes' theorem (\ref{thm:taubes}).} Let $\mathcal O_4$ be a four dimensional hyperbolic orbifold
whose topological fundamental group equals $G$;
it exists by Theorem~\ref{thm:orbifold}.
Let $V$ be the corresponding complex orbifold obtained from
$\mathcal O_4$ by the twistor construction.
Note that $\pi_1(V)\cong \pi_1(\mathcal O_4)$, since $V$ admits
a surjective map to $\mathcal O_4$ with connected and simply-connected fibers.

Finally, since $V$ is a complex analytic space we can resolve
its singularities by a theorem of Hironaka (for an expository account see \cite{Hau}),
then we apply Theorem~\ref{kollar}.
\qeds

\begin{wrapfigure}{r}{27mm}
\begin{lpic}[t(-5mm),b(3mm),r(0mm),l(0mm)]{pics/blows(0.5)}
\lbl[b]{24,41,;$\downarrow$ blow up of the $\downarrow$}
\lbl[t]{24,41,;blue curve}
\end{lpic}
\end{wrapfigure}

Alternatively, one can use the following explicit resolution.
This way one can avoid both the result of Hironaka and of Kollar.

Note, that the stabilizer of any orbi-point in $V$ in the above proof
is either $\mathbb Z_2$ or $\mathbb Z_2 \oplus \mathbb Z_2$.
Indeed, the action of the stabilizer on $\mathbb H^4$
preserves a complex structure on the tangent space
at the fixed point,
hence the stabilizer can not be $\mathbb Z_2\oplus \mathbb Z_2\oplus \mathbb Z_2$.
One can check that in appropriate local coordinates $(z_1,z_2,z_3)$ the action of $\mathbb Z_2$ is given by

$$(z_1,z_2,z_2)\mapsto (-z_1,-z_2,z_3),$$
while  the action of $\mathbb Z_2\oplus \mathbb Z_2$ is given by
$$(z_1,z_2,z_2)\mapsto (a_1\cdot z_1,a_2\cdot z_2,a_3\cdot z_3)$$
with $a_i=\pm 1$, $a_1\cdot a_2\cdot a_3=1$.

In particular the singularities with stabilizers $\mathbb Z_2\oplus \mathbb Z_2$ are formed by isolated points,
and at each of these points three complex curves with stabilizers $\mathbb Z_2$ meet
(each curve corresponds to a subgroup of order $2$ in $\mathbb Z_2\oplus \mathbb Z_2$).

Note that the singularity with stabilizer
$\mathbb Z_2 \oplus \mathbb Z_2$ can be represented locally by hypersurface
$w^2=x\cdot y\cdot z$ in $\mathbb C^4$.
Let us first blow up all irreducible components of the singular locus
that project to points on the hyperbolic orbifold $\mathcal O_4$.
(See the diagram.)
It is easy to see that
after this no singular component has a self-intersection, so we can consequently blow up each irreducible component.
This way we resolve all the singularities.

\parbf{Remark.} All complex three-folds that we construct are non-K\"ahler.
Indeed, a vertical fiber of the twistor space of $\HH^4$ has a non-compact family
of complex deformations with projections to $\HH^4$ of unbounded area. It follows that
for arbitrary Riemannian metric on $V$, complex deformations of a vertical fiber have unbounded
area as well. Hence neither $V$ nor any of its resolutions admit a K\"ahler metric.

\section{Comments}

\parbf{Optimizations of the construction.}
We had a lot of freedom in the above constructions;
as a rule we were choosing the way which is easier to write down.
Below we describe a few optimizations that may be used
elsewhere. In particular Lemma \ref{b3} is relevant for the work \cite{FP1} where
the existence compact symplectic Calabi--Yau six-manifolds with arbitrary fundamental
groups is deduced from Theorem \ref{thm:orbifold} of the present article.

\medskip

First, note that Proposition \ref{prop} holds even if $Y$ has only two orbi-points
of order $2$ in each cell;
this follows from Lemma~\ref{lem:2-points}.

\medskip

Second, in the proof of Proposition \ref{prop} one can relax partially condition 3. Namely,
it is sufficient that green and red curves cut $\Sigma$ into disks with holes.
Then in order to preform Step 2 of the proof one can appeal to the following:

\begin{thm}{Lemma}
 For any collection of positive integers $r_1,...,r_k$ there exists
a cover $\mathbb{S}^2\to \mathbb{S}^2$ of degree $n=\sum_i r_i$ ramified over $x_0,x_1,x_2,x_3\subset \mathbb{S}^2$
with ramifications of orders $r_i$ over $x_0$, and of orders at most two over $x_1,x_2,x_3$.
\end{thm}

This lemma holds since there exist
three involutions $\sigma_1,\sigma_2,\sigma_3$ in $S_n$ that
are composed altogether of $n+k-2$ transpositions, act transitively on
the set of $n$ elements, and such that $\sigma_1\cdot \sigma_2\cdot \sigma_3$ is a
product of disjoint cycles of lengths $a_i$.

Using the above remarks we can assume that the orbihedron $Y'$ constructed in
Proposition \ref{prop} has white
cells that are disks with arbitrary number $n$ of holes. To obtain such $Y'$
at Step 1 of the proof of Proposition \ref{prop} one can put in the interior of some white
cells a collection of $n$ disjoint couples of embedded red curves such that
the curves in each couple intersect in two points.

\medskip

Third, in the original construction, each two-cell of $Y'\subset \mathcal O_4$ is a topological disk.
It is easy to see that in $V$ (constructed in Section \ref{sec:taubes}) there are exactly two rational curves of singularities that project to each two-cell
(the interior of the cell has two preimages in both rational curves, and the boundary has one preimage).
On the other hand, if we modify the construction of
$Y'$ as above
then a cell in $Y'$ that is a disk with $n$ holes would correspond to two curves of genus $n$ in the corresponding $V$.
To summarize we have the following.

\begin{thm}{Lemma}\label{b3} For each integer $n>0$ and a finitely presented group $G$ there exists a compact oriented hyperbolic orbifold $\mathcal O_4$ with stabilizer $\ZZ_2^k$, $k=1,2,3$  having $\pi_1(\mathcal O_4)\cong G$, and such that the corresponding twistor space $V$ contains arbitrary large number of curves of $\ZZ_2$-singularities of genus $n$.
\end{thm}

\parbf{The construction in higher dimensions.}
Given a positive integer $m$,
consider the action $\Gamma_m\acts\HH^m$
defined by matrices with integer coefficients from $\QQ[\sqrt{5}]$
which preserve the quadratic form
$$\tfrac{1+\sqrt{5}}{2}\cdot x_0^2-x_1^2-\cdots-x_m^2.$$
The choice for $\QQ[\sqrt{5}]$ and $\tfrac{1+\sqrt{5}}{2}$ is made so that $\Gamma_2\acts\HH^2$ contains the Coxeter's action of right-angled regular pentagon.


We believe that the proof of Theorem~\ref{thm:orbifold} can be modified
to show that the action $\Gamma_m\z\acts\HH^m$ is telescopic
if the quotient space $\HH^m/\Gamma_m$
has finite fundamental group;
or, equivalently if a finite-index subgroup of $\Gamma_m$ is generated by elements of finite order.
According to \cite{bugaenko},
this holds at least for $m\le 7$;
in these dimensions $\Gamma_m$ contains a cocompact Coxeter's action.
Existence of a telescopic action on $\HH^6$
would lead via twistor construction \cite{FP}
to existence of symplectic Fano  orbifolds of dimension $12$ with arbitrary fundamental group.
We are not aware of any other applications of telescopic actions on $\HH^m$ for $m\ge 5$.

\end{document}